\documentclass[12pt,twoside]{article}

\usepackage{amsmath,amssymb,amsthm,hyperref}

\numberwithin{equation}{section}

\theoremstyle{plain}

\def\R{\mathbb{R}}
\def\N{\mathbb{N}}
\def\Nor{{\rm N}}

\def\balpha{{\boldsymbol\alpha}}
\def\bSigma{{\boldsymbol\Sigma}}
\def\diag{{\rm{diag}}}
\def\bzero{\boldsymbol{0}}

\def\bx{{\boldsymbol{x}}}

\def\d{\,{\rm{d}}}

\begin{document}

\title{\bf\Large A simple proof of the Gaussian correlation conjecture extended to multivariate gamma distributions}

\author{{\large T{\sc homas} R{\sc oyen}} \medskip \\
\centerline{University of Applied Sciences -- Bingen} \\
\centerline{Berlinstrasse 109, D-55411 Bingen, Germany} \medskip \\
\centerline{E-mail: \href{mailto:thomas.royen@t-online.de}{thomas.royen@t-online.de}}
}

\date{}
\maketitle

\begin{abstract}
An extension of the Gaussian correlation conjecture (GCC) is proved for multivariate gamma distributions (in the sense of Krishnamoorthy and Parthasarathy). The classical GCC for Gaussian probability measures is obtained by the special case with an integer degree of freedom $\nu = 1$.
\end{abstract}



\section{Introduction}
\label{sec:introduction}

Let $P$ be a probability measure on $\R^n$, $n > 1$, given by a Gaussian density,
$$
(2\pi)^{-n/2} |\bSigma|^{-1/2} \exp(-\tfrac12 \bx'\bSigma^{-1} \bx), \quad \bx \in \R^n,
$$
with a non-singular covariance matrix $\bSigma$. The Gaussian correlation conjecture (GCC) asserts the inequality 
\begin{equation}
\label{1.1}
P(C_1 \cap C_2) \ge P(C_1) P(C_2)
\end{equation}
for all convex and centrally symmetric sets $C_1, C_2 \subseteq \R^n$, (see \cite{dasgupta}). The bivariate case was proved in \cite{pitt}. Further milestones towards a complete proof were the papers \cite{schechtman} and \cite{harge}. In \cite{schechtman} the GCC is verified for all centered ellipsoids and for all sufficiently small $C_1, C_2$. In \cite{harge} a proof is given if only $C_1$ is a centered ellipsoid. Other recent claims to (very long) proofs can be found in \cite{memarian} and \cite{qing}. The proof presented here is totally different and comparatively short. 

According to \cite{schechtman} the GCC is equivalent to
\begin{equation}
\label{1.2}
P\Big(\bigcap_{i=1}^n A_i\Big) \ge P\Big(\bigcap_{i=1}^k A_i\Big) P\Big(\bigcap_{i=k+1}^n A_i\Big)
\end{equation}
with $A_i = \{|X_i| \le x_i\}$, $x_1,\ldots,x_n > 0$, $1 \le k < n$, and any $\Nor_n(\bzero,\bSigma)$-Gaussian vector $(X_1,\ldots,X_n)$. The inequality (1.2) was independently proved for $k=1$ in \cite{khatri} and \cite{schechtman}. Here, (1.2) is proved for $(X_1,\ldots,X_n)$ with an $n$-variate gamma distribution ($\Gamma_n(\alpha,R)$-distribution in the sense of Krishnamoorthy and Parthasarathy \cite{krishna}), defined by its Laplace transform 
\begin{equation}
\label{1.3}
|I_n+RT|^{-\alpha}
\end{equation}
with the identity matrix $I_n$, a non-singular correlation matrix $R = (\rho_{ij})$ and $T = \diag(t_1,\ldots,t_n)$, $t_1,\ldots,t_n \ge 0$.  Admissible values for $\alpha$ are $2\alpha \in \N$, all values $2\alpha > n-2$, $n \ge 2$, (for $n-2 < 2\alpha \le n-1$ see (1.7)) and all $\alpha > 0$ if the Laplace transform is infinitely divisible, for which sufficient and necessary conditions are found in \cite{bapat} and \cite{griffiths}. The $\Gamma_n(\alpha,R)$-distribution was originally derived from the joint distribution of the diagonal elements of a $W_n(2\alpha,R)$-Wishart matrix with the Laplace transform $|I_n+2RT|^{-\alpha}$. The classical GCC is obtained from the special case $\alpha = 1/2$. A proof of it could be easily extended inductively for increasing degrees of freedom $\nu = 2\alpha$ by means of convolution integrals. However, the proof provided here also includes non-integer values of $2\alpha$.

We need the non-central gamma probability density function (p.d.f.)
\begin{equation}
\label{1.4}
g_\alpha(x,y) = e^{-y} \sum_{k=0}^\infty g_{\alpha+k}(x) \frac{y^k}{k!}, \quad \alpha, x > 0, y \ge 0,
\end{equation}
with the central gamma p.d.f. $g_\alpha(x) = (\Gamma(\alpha))^{-1} x^{\alpha-1} e^{-x}$ and the corresponding non-central gamma cumulative distribution function (c.d.f.)
\begin{equation}
\label{1.5}
G_\alpha(x,y) = e^{-y} \sum_{k=0}^\infty G_{\alpha+k}(x) \frac{y^k}{k!}, \quad 
G_\alpha(x) = \int_0^x g_\alpha(\xi) \d \xi.
\end{equation}
If $R$ is represented by 
\begin{equation}
\label{1.6}
R = \lambda I_n + AA^T \iff \lambda^{-1} R = I_n + BB^T, \quad B = \lambda^{-1/2} A,
\end{equation}
with the minimal eigenvalue $\lambda$ of $R$ and an $n \times (n-1)$ matrix $A$ of rank $m \le n-1$, then the $\Gamma_n(\alpha,R)$-p.d.f. and the $\Gamma_n(\alpha,R)$-c.d.f., $2\alpha \in \N$ or $2\alpha > n-2$, can be represented by 
\begin{equation}
\label{1.7}
f(x_1,\ldots,x_n;\alpha,R) = E\Big(\prod_{j=1}^n \lambda^{-1} g_\alpha(\lambda^{-1} x_j,\tfrac12 b_jSb_j^T)\Big)
\end{equation}
and
\begin{equation}
\label{1.8}
F(x_1,\ldots,x_n;\alpha,R) = E\Big(\prod_{j=1}^n G_\alpha(\lambda^{-1} x_j,\tfrac12 b_jSb_j^T)\Big)
\end{equation}
respectively with the rows $b_j$ from $B$ and the expectation referring to $S$, a $W_{n-1}(2\alpha,I_{n-1})$-Wishart matrix.  These formulas are derived in a slightly more general form in \cite{royen95} and \cite{royen07}; they can be verified by a simple calculation using the Laplace transform of the non-central gamma distribution followed by integration over the Wishart (or pseudo-Wishart) distribution.

\section{Proof of the Gaussian correlation conjecture for \texorpdfstring{$\mathbf{\Gamma_n(\balpha,R)}$}--distributions}
\label{sec:proof}

\noindent{\large\bf Theorem 1.} \ \ 
{\it 
Let $R = \begin{pmatrix} R_{11} & R_{12} \\ R_{21} & R_{22}\end{pmatrix}$ be a non-singular partitioned $n \times n$-correlation matrix with $n_i \times n_i$ submatrices $R_{ii}$ and $R_{12}$ having positive rank.  Then, there holds for the $\Gamma_n(\alpha,R)$-c.d.f.
\begin{equation}
\label{2.1}
F(x_1,\ldots,x_n;\alpha,R) > F(x_1,\ldots,x_{n_1};\alpha,R_{11}) F(x_{n_1+1},\ldots,x_n;\alpha,R_{22})
\end{equation}
for all positive numbers $x_1,\ldots,x_n$ and $2\alpha \in \N$ or $2\alpha > n-2$.
}

\bigskip

\noindent
{\bf Proof.}  All the matrices
\begin{equation}
\label{2.2}
R_\tau = \begin{pmatrix} R_{11} & \tau R_{12} \\ \tau R_{21} & R_{22}\end{pmatrix}, \quad 0 \le \tau \le 1,
\end{equation}
are non-singular correlation matrices. Theorem 1 will be proved if we show 
\begin{equation}
\label{2.3}
\frac{\partial}{\partial\tau} F(x_1,\ldots,x_n;\alpha,R_\tau) > 0, \quad 0 < \tau < 1.
\end{equation}
By means of the Laplace transform we shall represent the left-hand side of (2.3) by a sum of not identically vanishing non-negative terms.

Let $A_J$ denote the submatrices with row and column indices $i\in J \subseteq \{1,\ldots,n\}$ from any $n \times n$-matrix $A$. The determinant $|I_n+RT|$ is equal to $1 + \sum_J |R_J| |T_J|$, $\emptyset \neq J \subseteq \{1,\ldots,n\}$.  With $J = J_1 \cup J_2$, $J_1 = J \cap \{1,\ldots,n_1\} \neq \emptyset$, $J_2 = J \cap \{n_1+1,\ldots,n\} \neq \emptyset$, $R_{\tau,J} = \begin{pmatrix} R_{J_1} & \tau R_{J_1,J_2} \\ \tau R_{J_2,J_1} & R_{J_2}\end{pmatrix}$, $r_{J_1,J_2} = {\rm{rank}}(R_{J_1,J_2})$, and the squared canonical correlations $\rho_{J_1,J_2,i}^2$, $i=1,\ldots,r_{J_1,J_2}$, which are the positive eigenvalues of $R_{J_1}^{-1/2} R_{J_1,J_2} R_{J_2}^{-1} R_{J_2,J_1} R_{J_1}^{-1/2}$, we find 
\begin{align*}
|R_{\tau,J}| &= |R_{J_1}| \; |R_{J_2}| \; |I_{J_1} - \tau^2 R_{J_1}^{-1/2} R_{J_1,J_2} R_{J_2}^{-1} R_{J_2,J_1} R_{J_1}^{-1/2}| \\
&= |R_{J_1}| \; |R_{J_2}| \; \prod_{i=1}^{r_{J_1,J_2}}(1-\tau^2\rho_{J_1,J_2,i}^2),
\end{align*}
and 
\begin{align*}
\frac{\partial}{\partial\tau} |R_{\tau,J}| 
&= |R_{J_1}| \; |R_{J_2}| \; \frac{\partial}{\partial\tau} \prod_{i=1}^{r_{J_1,J_2}}(1-\tau^2\rho_{J_1,J_2,i}^2) \\
&= -2\tau |R_{\tau,J}| \; \sum_{i=1}^{r_{J_1,J_2}} \frac{\rho_{J_1,J_2,i}^2}{1-\tau^2\rho_{J_1,J_2,i}^2}.
\end{align*}
Therefore, 
\begin{align}
\label{2.4}
\frac{\partial}{\partial\tau} |I_n+R_\tau T|^{-\alpha} &= 2\alpha\tau |I_n+R_\tau T|^{-(\alpha+1)} \sum_J |R_{\tau,J}| \Big(\sum_{i=1}^{r_{J_1,J_2}} \frac{\rho_{J_1,J_2,i}^2}{1-\tau^2\rho_{J_1,J_2,i}^2}\Big) \prod_{j \in J} t_j \nonumber \\
&= |I_n+R_\tau T|^{-(\alpha+1)} \sum_J c_J(\tau) \prod_{j \in J} t_j,
\end{align}
which is the Laplace transform $h^*(t_1,\ldots,t_n;\alpha,R_\tau)$ of 
\begin{equation}
\label{2.5}
h(x_1,\ldots,x_n;\alpha,R_\tau) = \sum_J c_J(\tau) \Big(\prod_{j \in J} \frac{\partial}{\partial x_j}\Big) f(x_1,\ldots,x_n;\alpha+1,R_\tau),
\end{equation}
where 
$$
c_J(\tau) = 2\alpha\tau \; |R_{\tau,J}| \sum_{i=1}^{r_{J_1,J_2}} \frac{\rho_{J_1,J_2,i}^2}{1-\tau^2\rho_{J_1,J_2,i}^2}.
$$

We have to verify that 
\begin{equation}
\label{2.6}
h(x_1,\ldots,x_n;\alpha,R_\tau) = \frac{\partial}{\partial\tau} f(x_1,\ldots,x_n;\alpha,R_\tau).
\end{equation}
Integration of the left-hand side in (2.4) over $[0,\tau]$ provides 
\begin{align}
\label{2.7}
|I_n+&R_\tau T|^{-\alpha} - |I_n+R_0 T|^{-\alpha} \nonumber \\
&= \int_{\R_+^n} \big(f(x_1,\ldots,x_n;\alpha,R_\tau) - f(x_1,\ldots,x_n;\alpha,R_0)\big) \prod_{j=1}^n e^{-t_j x_j} \d x_j \nonumber \\
&= \int_0^\tau h^*(t_1,\ldots,t_n;\alpha,R_\vartheta) \d\vartheta \nonumber \\
& = \int_{\R_+^n} \Big(\int_0^\tau h(x_1,\ldots,x_n;\alpha,R_\vartheta) \d\vartheta \Big) \prod_{j=1}^n e^{-t_j x_j} \d x_j,
\end{align}
where the change of the order of integration can be justified by Fubini's criterion in the following way:

Applying (1.7), the identity
$$
\frac{\partial}{\partial x}g_{\alpha+1}(x,y) = g_\alpha(x,y) - g_{\alpha+1}(x,y),
$$
decompositions $\lambda_\vartheta^{-1} R_\vartheta = I_n + B_\vartheta B_\vartheta^T$ - similar to (1.6) - with rows $b_{j\vartheta}$ in $B_\vartheta$, indicator functions $e_J$ of $J$,  and the cardinalities $|J|$ of $J$, we obtain 
\begin{align*}
|h(&x_1,\ldots,x_n;\alpha,R_\vartheta)| \\
& \le \sum_J c_J(\vartheta) \Big|\Big(\prod_{j \in J} \frac{\partial}{\partial x_j}\Big) f(x_1,\ldots,x_n;\alpha+1,R_\vartheta)\Big| \\
& \le \sum_J c_J(\vartheta) \lambda_\vartheta^{-|J|} E\Bigg(\prod_{j=1}^n \lambda_\vartheta^{-1} \Big(e_J(j) \; g_{\alpha}\big(\lambda_\vartheta^{-1} x_j,\tfrac12 b_{j\vartheta} Sb_{j\vartheta}^T\big) + g_{\alpha+1}\big(\lambda_\vartheta^{-1} x_j,\tfrac12 b_{j\vartheta} S b_{j\vartheta}^T\big)\Big)\Bigg),
\end{align*}
$S \sim W_{n-1}(2(\alpha+1),I_{n-1})$.  

Therefore, the Laplace transform of $|h|$ is bounded by a linear combination of integrals of the form
$$
\int_{\R_+^n} E\Big(\prod_{j=1}^n \lambda_\vartheta^{-1} g_{\alpha+e_K(j)}\big(\lambda_\vartheta^{-1} x_j,\tfrac12 b_{j\vartheta} Sb_{j\vartheta}^T\big)\Big) \prod_{j=1}^n e^{-t_j x_j} \d x_j = |I_n + R_\vartheta T|^{-(\alpha+1)} \prod_{i\in K^c} (1+\lambda_\vartheta t_i),
$$
$K \subseteq \{1,\ldots,n\}$, where the right-hand side of this equation is obtained in a way similar to the Laplace transform of $f(x_1,\ldots,x_n;\alpha,R)$ in (1.7).  The coefficients of this linear combination are non-negative continuous functions of $\vartheta \in [0,\tau]$. Integration over $[0,\tau]$ yields a finite value, which implies (2.7) and (2.6). 

Finally, integration over $x_1,\ldots,x_n$ in (2.5) leads to 
$$
\frac{\partial}{\partial\tau} F(x_1,\ldots,x_n;\alpha,R_\tau) = \sum_J c_J(\tau) \Big(\prod_{j \in J} \frac{\partial}{\partial x_j} \Big) F(x_1,\ldots,x_n;\alpha+1,R_\tau) > 0,
$$
since the canonical correlations in (2.4) do not identically vanish because of the positive rank of $R_{12}$. $\qed$

\medskip

\noindent
{\bf Remarks.}  By continuity, Theorem 1 holds also for a singular $R$, at least with ``$\ge$'' instead of ``$>$''.  

A modified proof uses the characteristic function $\hat{h}$ corresponding to (2.4). It is the continuous limit of the characteristic functions of the signed measures defined by $\epsilon^{-1}\big(F(x_1,\ldots,x_n;\alpha,R_{\tau+\epsilon}) - F(x_1,\ldots,x_n;\alpha,R_\tau)\big)$, $\epsilon \to 0$. Thus, $\hat{h}$ is the characteristic function of the signed measure which is defined by $\frac{\partial}{\partial \tau}F(x_1,\ldots,x_n;\alpha,R_\tau)$ and which is the difference of two finite positive measures on $\R^n$.  

An approximation by a series of univariate integrals for the difference between the two sides in (2.1) is proposed in \cite{royen13} with identical values $x_i = x$ under the additional conditions
\begin{align*}
\rho_1 &= \frac{2}{n_1(n_1 - 1)} \sum_{1 \le i < j \le n_1} \rho_{ij}> 0, \quad
\rho_2 = \frac{2}{n_2(n_2 - 1)} \sum_{n_1+1 \le i < j \le n} \rho_{ij} > 0, \\
\rho^2 &= \Big(\frac{1}{n_1n_2}\sum_{i=1}^{n_1}\sum_{j=n_1+1}^n \rho_{ij}\Big)^2 \le \rho_1 \rho_2.
\end{align*}
The error of this approximation tends to zero with a decreasing variability of the correlations within the submatrices $R_{11}$, $R_{22}$, and $R_{12}$. This approximation is recommended in particular for small values of $1-F(x,\ldots,x;\alpha,R)$. In a modified
approximation, $\rho^2$ is replaced by $\frac{1}{n_1n_2}\sum_{i=1}^{n_1}\sum_{j=n_1+1}^n \rho_{ij}^2$.

\vskip0.5truein

\noindent
{\bf Acknowledgement.} The author wishes to thank Donald Richards from the Pennsylvania State University for several valuable recommendations which made this paper easier to read.

\vskip 0.5truein

\bibliographystyle{ims}

\begin{thebibliography}{99}

\bibitem{bapat}
Bapat, R.B. (1989). Infinite divisibility of multivariate gamma distributions and M-matrices, {\it Sankhy$\bar{\rm{a}}$} 51, 73--78.

\bibitem{dasgupta}
Das Gupta, S., Eaton, M. L., Olkin, I., Perlman, M., Savage, L.J., and Sobel, M. (1972).  Inequalities on the probability content of convex regions for elliptically contoured distributions.  {\sl Proc. Sixth Berkeley Symp. Math. Statist. Probab.} 2, 241--264, Univ. of California Press, Berkeley.

\bibitem{griffiths}
Griffiths, R.C. (1984). Characterization of infinitely divisible multivariate gamma distributions, {\it J. Multivariate Anal.} 15, 13--20.

\bibitem{harge}
Harg\'e, G. (1999).  A particular case of correlation inequality for the Gaussian measure.  {\it Ann. Probab.} { 27}, 1939--1951.

\bibitem{khatri}
Khatri, C.G. (1967). On certain inequalities for normal distributions and their application to simultaneous confidence bounds, {\it Ann. Math. Stat.} 38, 1853--1867.

\bibitem{krishna}
Krishnamoorthy, A.S., and Parthasarathy, M. (1951).  A multivariate gamma-type distribution.  {\it  Ann. Math. Stat.} { 22}, 549--557.

\bibitem{memarian}
Memarian, Y. (2013). The Gaussian correlation conjecture proof, arXiv:1310. 8099v1 [mathPR].

\bibitem{pitt}
Pitt, L.D. (1977).  A Gaussian correlation inequality for symmetric convex sets. {\it Ann. Probab.} { 5}, 470--474.

\bibitem{qing}
Qingyang, G. (2013). The Gaussian correlation inequality for symmetric convex sets, arXiv:1012.0676v4 [mathPR].

\bibitem{royen95}
Royen, T. (1995). On some central and non-central multivariate chi-square distributions, {\it Statist. Sinica} 5, 373--397.

\bibitem{royen07}
Royen, T. (2007). Integral representations and approximations for multivariate gamma distributions, {\it Ann. Inst. Statist. Math.} 59, 499--513.

\bibitem{royen13}
Royen, T. (2013). Some upper tail approximations for the distribution of the maximum of correlated chi-square or gamma random variables, {\it Far East J. Theor. Statist.} 43, 27--56.

\bibitem{schechtman}
Schechtman, G., Schlumprecht, T., and Zinn, J. (1998).  On the Gaussian measure of the intersection, {\it Ann. Probab.} { 26}, 346--357.

\bibitem{sidak}
\v Sid\'ak, Z. (1967).  Rectangular confidence regions for the means of multivariate normal distributions, {\it J. Amer. Statist. Assoc.} 62, 626--633.
\end{thebibliography}

\end{document}